\documentclass[12pt]{amsart}
\usepackage{amsmath,amssymb,epsfig,amscd,xy}
\xyoption{all}

\newtheorem{theorem}{Theorem}

\newtheorem{definition}[theorem]{Definition}

\newtheorem*{main}{Theorem}
\newtheorem*{main-lemma}{Lemma}

\newcommand{\cE}{{\mathcal E}}

\newcommand{\bG}{{\mathbb G}}

\DeclareMathOperator{\Supp}{Supp}
\DeclareMathOperator{\SL}{SL}
\DeclareMathOperator{\Spec}{Spec}
\DeclareMathOperator{\Norm}{Norm}
\DeclareMathOperator{\N}{N}
\DeclareMathOperator{\Gal}{Gal}

\newcommand{\bP}{{\mathbb P}}

\newcommand{\bC}{{\mathbb C}}

\newcommand{\bQ}{{\mathbb Q}}
\newcommand{\lra}{\longrightarrow}

\newfont{\cyrr}{wncyr10}

\begin{document}

\title{ One half log discriminant } \author{Lucien Szpiro and Thomas
  J.~Tucker}

\keywords{Equidistribution, Mahler measure, elliptic curves} 
\subjclass[2000]{Primary 14G40, Secondary 11G50, 11G05}

\address{
Lucien Szpiro \\
Ph.D. Program in Mathematics\\
Graduate Center of CUNY\\
365 Fifth Avenue\\
New York, NY 10016-4309
}

\email{lszpiro@gc.cuny.edu}

\address{
Thomas Tucker\\
Department of Mathematics\\
Hylan Building\\
University of Rochester\\
Rochester, NY 14627
}

\email{ttucker@math.rochester.edu}

\thanks{The first author was partially supported
  by NSF Grant 0071921.  The second author was partially supported by
  NSF Grant 0101636}

\begin{abstract}
  We give a geometric proof that one may compute a particular
  generalized Mahler integral using equidistribution of preperiodic
  points of a dynamical system on the sphere. The dynamical system is
  associated to the multiplication by 2 map on an elliptic curve over
  a number field $K$ with Weierstrass equation $y^2 = P(x)$ (a
  Latt{\`e}s dynamical system).  At each finite place $v$, we prove
  the local equidistribution formula
  $$v(\Delta) =\lim_{n\rightarrow\infty}\frac{1}{n^2} (D.H_n)_v,
  $$
  where $H_n$ is the Zariski closure in
  $\mathbb{P}^1_{\mathcal{O}_K}$ of the image in $\bP^1_K$ of the
  $n$-torsion minus the 2-torsion and $\Delta$ is the discriminant of
  the polynomial $P(x)$.  One consequence of this result is the
  formula
 $$\frac{1}{2} \log|\Norm_{K/\bQ} (\triangle)| =\sum_{\sigma:
   K\rightarrow\mathbb{C}} \int_{\mathbb{P}^1(\mathbb{C})} \log
 |P(x)|_\sigma \hspace{1mm} \frac{dx \wedge d\bar{x}}{\Im(\tau)_\sigma
   |P(x)|_\sigma^2 }.$$
\end{abstract}

\maketitle

In \cite{SUZ}, Szpiro, Ullmo, and Zhang proved that for any abelian
variety $A$ over $\bQ$, any continuous function $g$ on $A(\bC)$, and
any nonrepeating sequence of point $\beta_n \in A({\overline \bQ})$ with
N{\'e}ron-Tate height tending to zero, one has
$$
\frac{1}{|\Gal(\beta_n)|}\sum\limits_{\sigma \in{\Gal(\beta_n)}}
g(\sigma(\beta_n)) = \int_{A(\bC)} g \, d \mu,$$ where $d \mu$ is
the normalized Haar measure on $A$ and $\Gal(\beta_n)$ is the
Galois group of the Galois closure of $\bQ(\beta_n)$ in $\bC$.
This result says, in effect, that Galois orbits of points with
small N{\'e}ron-Tate height are equidistributed in $A(\bC)$. Ullmo
\cite{ullmobogomolov} and Zhang \cite{zhangbogomolov} later used
this fact to give proofs of the Bogomolov conjecture for abelian
varieties.

When the abelian variety $A$ is an elliptic curve, the
multiplication by 2 map gives rise to a map on the projective line,
called a Latt{\`e}s map.  Thus, in this case, the work of \cite{SUZ}
can be viewed as an equidistribution result for a rational map on the
projective line.  Recently, a variety of authors have proven more
general equidistribution results for arbitrary rational maps of degree
greater than 1 on the projective line; see Autissier \cite{autissier},
Baker/Rumely \cite{BR}, Bilu \cite{bilu}, Chambert-Loir \cite{CL}, and
Favre/Rivera-Letelier \cite{FR1,FR2}, for example.  Many of these
results hold for measures at finite places as well as at archimedean
places.

In \cite{szmalher}, it is shown that for any nonconstant map
$\varphi: \bP^1 \lra \bP^1$ of degree greater than 1 over a number
field $K$, the canonical height $h_\varphi(\alpha)$ of an
algebraic point $\alpha$ with minimal polynomial $F$ can be
calculated by integrating $\log |F|$ along the invariant measures
for the map $\varphi$.  This gives a generalization of the notion
of a Mahler measure of a polynomial (see \cite{mahler}).  Everest,
Ward, and Fhlathuin \cite{algebraicdynamics, elliptic} had
previously extended the notion of Mahler measures to elliptic
curves.

Additional difficulties arise, however, when one attempts to prove
equidistribution results for the functions $\log |F|_v$.  Indeed,
the exact analog of the main result of \cite{SUZ} is not true when
the continuous function $g$ is replaced by a function of the form
$\log |F|$ (see \cite{aut-letter} or \cite{BIR}).  On the other
hand, it is possible to prove an equidistribution result for
functions of the form $\log |F|_v$ provided that one averages over
all points of period $n$ as $n$ goes to infinity rather than over
Galois orbits of families of points of small height (see
\cite{sztuequi}). In the case of elliptic curves, Baker, Ih, and
Rumely \cite{BIR} were able to refine this to prove that for any
algebraic number $\alpha$ and any Latt{\`e}s map $\varphi$ one has
$$
[K(\alpha):\bQ] \,h_\varphi(\alpha) = \sum_{\text {places $v$ of
    $K$}} \lim_{n \to \infty} \frac{1}{|\Gal(\beta_n)|}
\sum\limits_{\sigma \in{\Gal(\beta_n)}} \log |F(\beta_n^\sigma)|_v
$$
for any nonrepeating sequence of algebraic points $\beta_n$ such
that $h_\varphi(\beta_n) = 0$ for all $n$.  Both \cite{BIR} and
\cite{sztuequi} use results from diophantine approximation,
specifically Roth's theorem (\cite{Roth}) and A.~Baker's work on
linear forms in logarithms (\cite{Baker}).

When one applies the results of \cite{BIR} and \cite{szmalher} to the
points of period 2 for a Latt{\`e}s map corresponding to
multiplication by 2 on the elliptic curve $E$ with Weierstrass
equation $y^2= P(x)$, one obtains the formula
\begin{equation*}
\begin{split}
  \frac{1}{2} \log|\Norm_{K/\bQ} (\triangle)| & = \sum_{\sigma: K
    \hookrightarrow \mathbb{C}} \lim_ {n\rightarrow\infty}\frac{1}
  {n^2} \log \prod_{\beta \in \Supp H_n} |P(\beta)|_\sigma\\
  & = \sum_{\sigma: K \hookrightarrow\mathbb{C}}
  \int_{\mathbb{P}^1(\mathbb{C})} \log |P(x)|_\sigma \hspace{1mm}
  \frac{dx \wedge d\bar{x}}{\Im(\tau)_\sigma
    |P(x)|_\sigma^2 }, \\
\end{split}
\end{equation*}
where $\Delta$ is the discriminant of $F$ over $K$ and $\tau_\sigma$
denotes the element corresponding to the elliptic curve $E_\sigma$ in
the fundamental domain for the action of $\SL(2,\mathbb{Z})$ on the
Poincare upper half space in $\mathbb{C}$.  Using the product formula
and the fact that $h_\varphi(\alpha) = 0$ for periodic points
$\alpha$, this is equivalent to showing that at each nonarchimedean
place $v$ of $K$, we have
\begin{equation*}
\begin{split}
  \lim_ {n\rightarrow\infty}\frac{1} {n^2} \log
  \prod_{\beta \in \Supp H_n} |P(\beta)|_\sigma &
  = \lim_{n\rightarrow\infty}\frac{1}{n^2}\log
  \Norm_{H_n/\mathcal{O}_K}(P_{|H_n}) \\
  & = \lim_{n\rightarrow\infty}\frac{1}{n^2} \sum_v (D.H_n)_v \log
  \N(v),
\end{split}
\end{equation*}
where $N(v)$ is the cardinality of the residue field at $v$ and $H_n$
is the Zariski closure in $\mathbb{P}^1_{\mathcal{O}_K}$ of the image
in $\bP^1_K$ of the $n$-torsion minus the 2-torsion

We give here a a local proof using blow-ups of closed points and
intersection theory on $\mathbb{P}^1_V$. This proof uses
resolution of singularities (in fact separation of branches) of
one-dimensional schemes by blowing up. It also uses information
about the special fiber of an elliptic curve with semistable
reduction (see \cite{silveradvanced}).  We do not use diophantine
approximation.  The proof we give is valid for equicharacteristic
$V$ (local geometric case) as well as for unequal characteristic
(local arithmetic case). Note that the case of positive
characteristic does not follow from the results of \cite{BIR} and
\cite{sztuequi}, since the relevant approximation theorems are not
valid in characteristic $p$.  Relations between the discriminant
of an elliptic curve and its $n$-torsion have been studied in
\cite{szgroth} and in \cite {szpiroullmo}.  The main theorem of
this paper is the following.

\begin{main}
  Let $V$ be a discrete valuation ring with fraction field $K$.  Let
  $y^2 = P(x)$ be the minimal Weierstrass equation with coefficients
  in $V$ of an elliptic curve $E$. Suppose that $E$ has semi-stable
  reduction over $V$. Let $D$ be the scheme of zeroes of $P(x)$ in
  $\mathbb{P}^1_V$ and let $H_n$ be the Zariski closure in
  $\mathbb{P}^1_V$ of the image in $\mathbb{P}^1_K$ of the kernel of
  the multiplication by $n$ in $E_K$ minus the 2-torsion.  Then
  $$\lim_{n \rightarrow\infty} \frac{1}{n^2}(D.H_n)_v =
  \frac{1}{2}v(\triangle),$$
  where $v$ is the normalized valuation of
  $V$, $\triangle$ is the discriminant of $E$ over $V$, and $(-.-)_v$
  is the geometric intersection pairing on the surface
  $\mathbb{P}^1_V$.
\end{main}

For simplicity we we will assume that the roots of $P(x)$ are rational
over $K$. Also for simplicity we assume that the residual
characteristic of $V$ is not 2. These two conditions are not essential
for the theorem but they make the proof easier.  The valuation of the
discriminant is then even; we write $2k = v(\triangle)$. One knows
(see for example \cite{silveradvanced} or \cite{deligneformulaire})
that the closed fiber of the minimal model $\cE$ for $E$ over $V$ is a
cycle of $2k$ projective lines of self intersection $(-2)$; it is
obtained by blowing up the plane model for the elliptic curve $k$
times.  Recall that the N{\'e}ron model in this case is
$$ \cE \setminus \{ \text {singular points of the special
fiber} \} \text {  (see figure 3)}.$$

The strategy of the proof is to compute intersections in the N{\'e}ron
model for $E$ after a suitable base change.  The multiplicative
structure of the special fiber is simply $\bG_m$ crossed with the
group of components.  One can easily see how the $n$-torsion
distributes itself among the components, and that allows one to
calculate intersections without difficulty.

We will use the fact that the N{\'e}ron model for $E$ has $2k$
components in its special fiber(see \cite[page 365]{silveradvanced}).
It is naturally a $2:1$ cover of a model for $\bP^1$ with $k+1$
components. The hyperelliptic map induces a map on components that
sends inverse component and its inverse to a single component; there
are two components that are their own inverse (the identity and the
component of order 2), which gives a total of $(k-1)/2 + 2 = k+1$
components on a model of $\bP^1$.

We begin with the plane model $E$ for $E_K$ coming from the
equation $y^2 = P(x)$.
\\[0.3 in]
\includegraphics[scale=0.5]{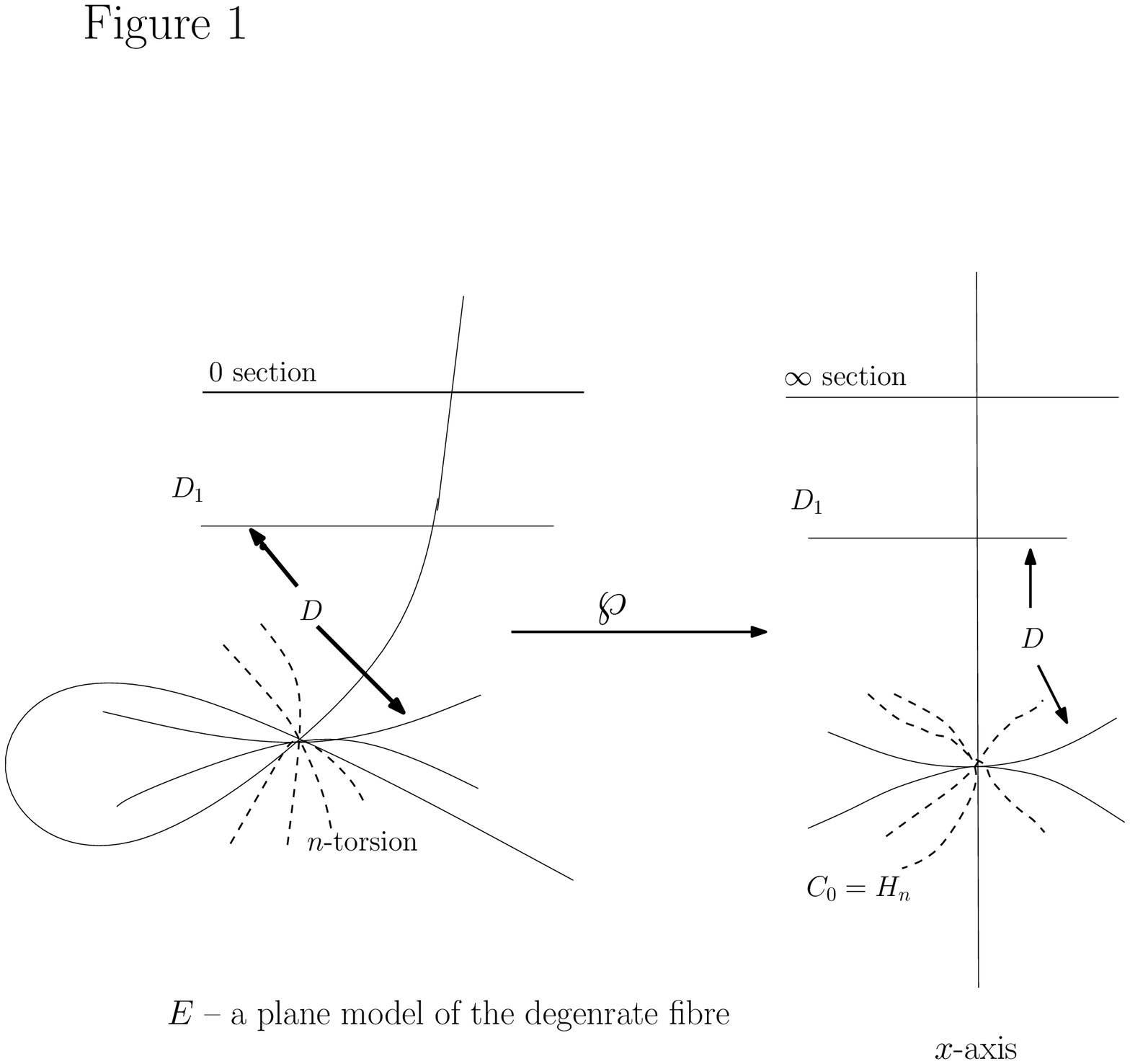}
\\[0.3 in]

   \begin{definition}  Let $D_0$ denote the divisor
     $D$.  We define the divisor $D_i$ recursively (for $i \leq k$) as
     the proper transform of $D_{i-1}$ for the blow-up $\sigma_i: X_i
     \rightarrow X_{i-1}$ centered at the point $P_{i-1}$ of
     multiplicity $2$ on $D_{i-1}$.
\end{definition}

Note that this is a horizontal divisor of degree 3 intersects the
special fiber $F_0$ of $\mathbb{P}^1_V = X_0$ in 2 points: one $P_0$
of multiplicity 2 on $D_0$, the other one of multiplicity 1 on $D_0$.
We now define the divisors in our models $X_i$ coming from $H_n$.

\begin{definition}The horizontal divisor $C_0$
  is defined to be $H_n$ for some fixed odd $n$. The divisor $C_i$ is
  the proper transform of $C_{i-1}$ in $X_i$.
\end{definition}

The degree of $C_0$ is $(n^2-1)/2$ when $n$ is odd and $(n^2/2) - 3$
when $n$ is even.  This follows from the fact that the hyperelliptic
map sends each point and its inverse to the same point in $\bP^1$.

\begin{definition}Define $\wp_K: E_K \rightarrow \mathbb{P}^1_K $ to be the projection
  onto ``the $x$ axis" (i.e., $\wp$ is the Weierstrass $\wp$
  function).  We will, by abuse of language, note $\wp: E_i\rightarrow
  X_i$ to be the the extension of $\wp_K$ to model $E_i$ for $E_K$
  over $V$.
\end{definition}

The figure 2 illustrates the situation.
\\[0.2 in]
\includegraphics[scale=0.6]{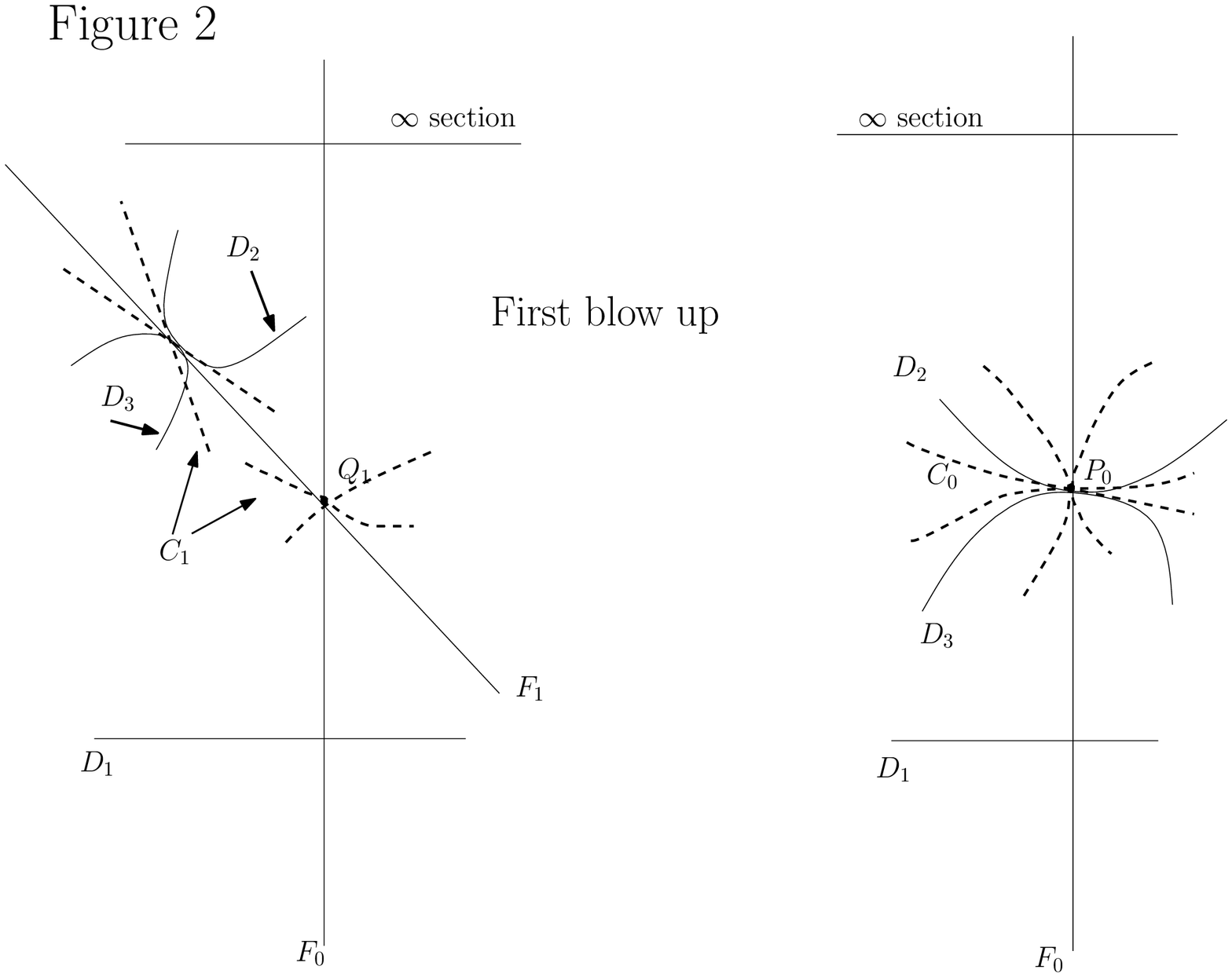}
\\[0.2 in]

\begin{main-lemma}
  Assume that $n$ is odd or that the residual characteristic is not 2,
  then after $k$ successive blow-ups of the points $P_i$ of
  multiplicity 2 on $D_i$, the proper transform $D_k$ is {\'e}tale and
  the proper transforms $D_k$ and $C_k$ do not meet.
\end {main-lemma}
\begin{proof}({\it Of Lemma.}) If $\wp^*(C_k)$ and
  $\wp^*(D_k)$ are both in the N{\'e}ron model (i.e., if $n$ and
  $(2k)$ have a common factor $m$), then $H_n$ and 2-torsion are
  distinct; hence, when the characteristic is not 2, they do not meet
  in the N{\'e}ron model.  If $n$ is prime to $2k$ and $\wp^*(H_k)$ is
  not inside the N{\'e}ron model, then $\wp^*(H_k)\bigcap \wp^*(D_k) =
  \emptyset$, since $\wp^*(D_k)$ is in the N{\'e}ron model (see figure
  3).
\end{proof}
\includegraphics[scale=0.6]{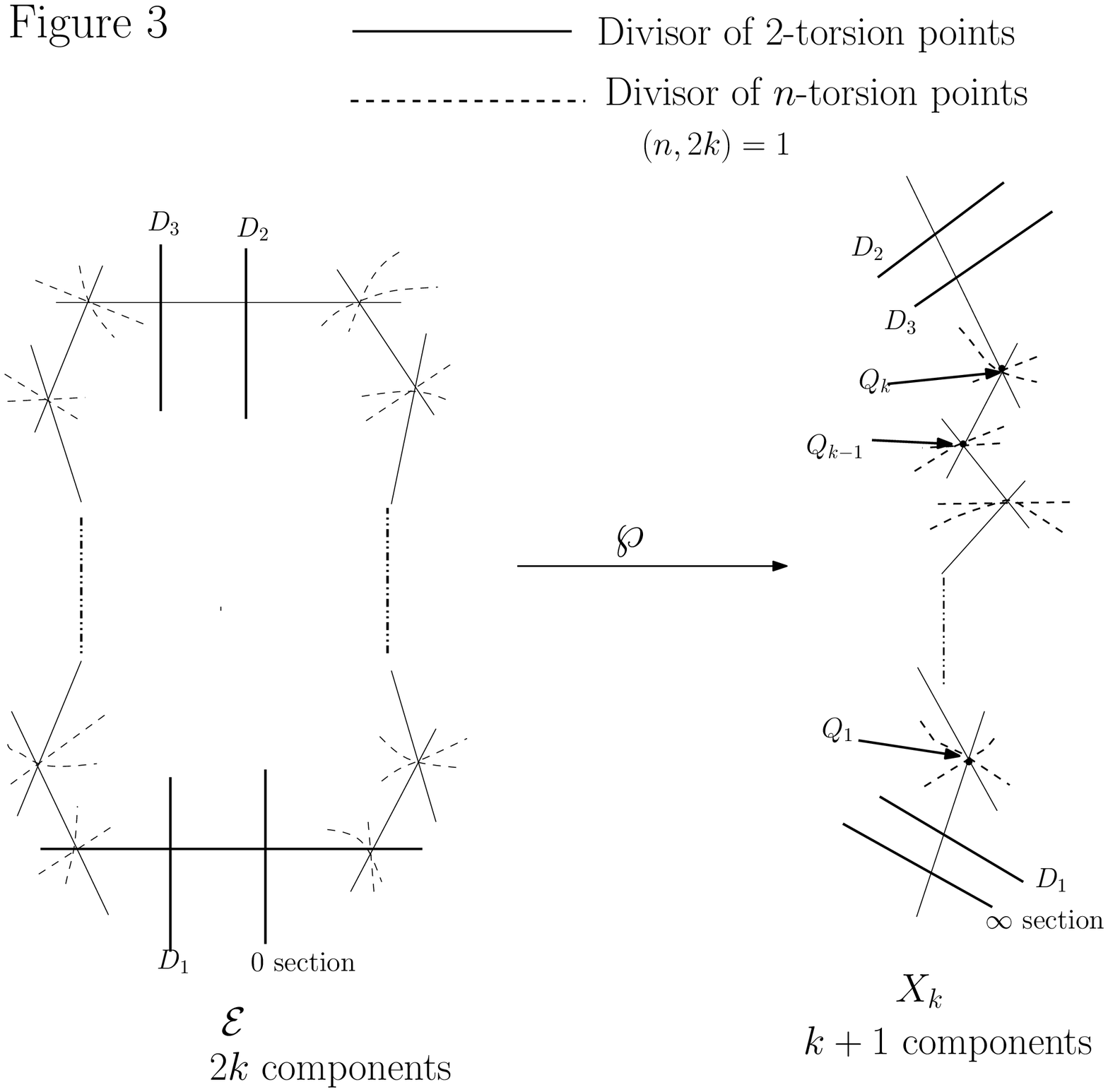}
\\[0.4 in]

We are now ready to prove the main theorem.
\begin{proof}
  We treat first the case when $n$ is prime to $2k$.  The exceptional
  divisor of $\sigma_i$ will be denoted as $F_i$. By abuse of language the
  proper transform of $F_i$ will still be called $F_i$ after
  $\sigma_{i+1},\dots,\sigma_k$. We will let $Q_i$ denote the point of
  intersection of $F_i$ with $F_{i-1}$ in $X_i$ (see figure 2).

  We will denote the usual pull-back map for divisors with $*$. We
  denote the composed map $\sigma_i \cdot \sigma_{i-1} \cdots
  \sigma_1$ as $\rho_i$.  After $i$ blow-ups, one has integers
  $m_{j,i}$ such that
$$\rho_i^*D=D_i+ \sum_{j\leq i}m_{j,i}F_j$$
  and
  $$\pi_{i+1}^*D =\sigma_{i+1}^* D_i+ \sum_{j\leq i} m_{j,i}
  \sigma_{i+1}^*F_j.$$
\\[0.4 in]
  \includegraphics[scale=0.6]{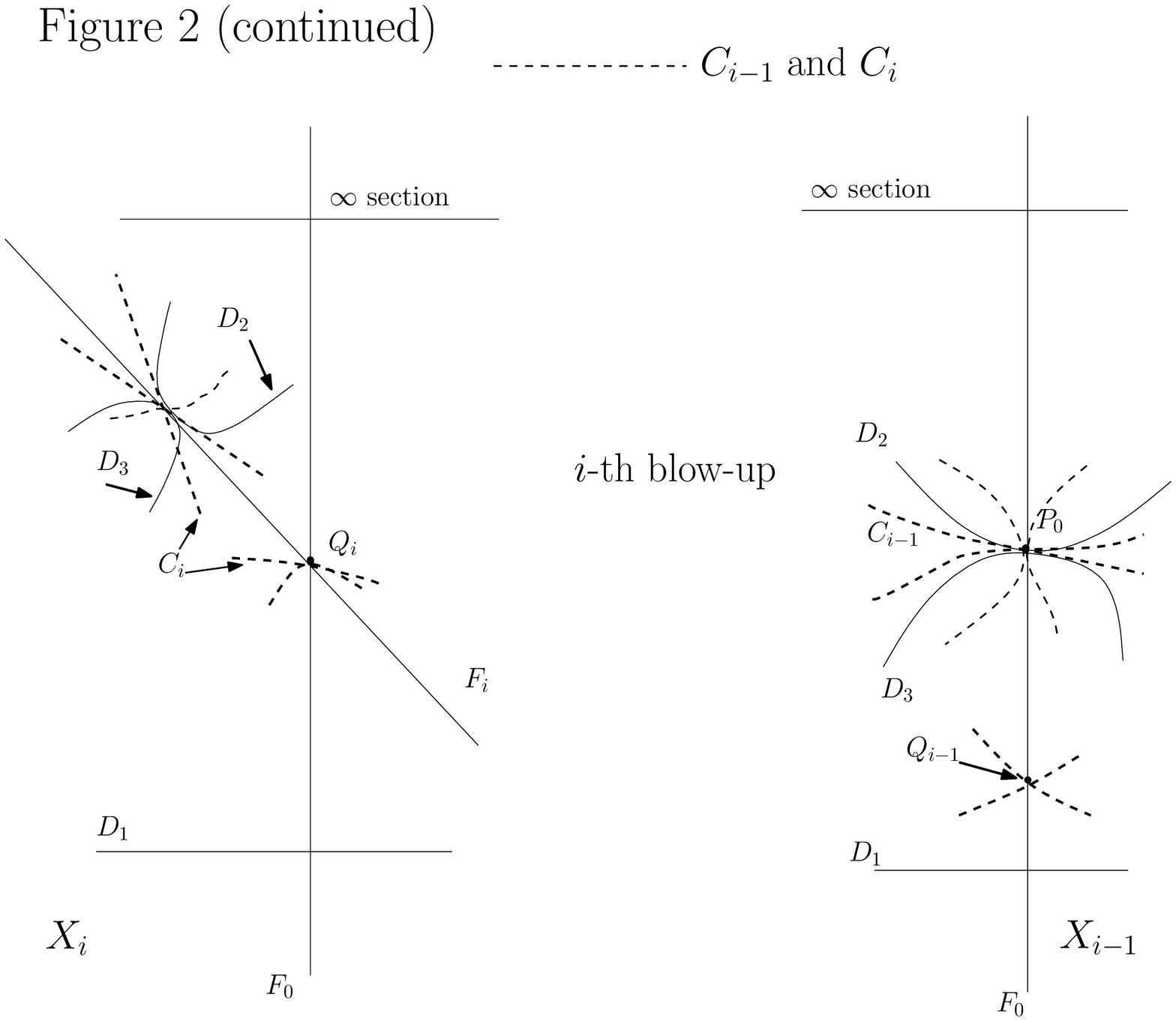}
\\[0.4 in]

  As long as $i$ is less than $(k-1)$, one has
  $$\sigma_{i+1}^*D_i=D_{i+1}+2F_{i+1},$$
  since the multiplicity of
  $D_i$ at $P_i$ is still 2.  Since $\sigma_{i+1}^*F_i=F_i+F_{i+1}$ we
  have $m_{j,i}=m_{j,j}$ for any $i \geq j$, so we have
  $$m_{j, i} = m_{j,j}= m_{j - 1, j - 1}+2$$
  for all $i \geq j$.
  Thus, by induction, we have $m_{j,j} = 2j$ for each $j$, which means
  that $m_{j,i}=2j$ for all $i \geq j$.

  The intersection multiplicity we are looking for can be computed as
  follows
  $$(*)(H_n.D) = (C_0.D_0)=(C_k.(D_k+\sum_{j\leq k}m_{j,k}F_j))=
  2\sum_{j\leq k} j(C_k.F_j).$$
  One is left with computing each
  $(C_k.F_j)$. We will achieve this by looking at the special fiber of
  various models of $E_K$ over $V$. By the projection formula for
  $\wp$ we can compute intersections on the minimal model $\cE$ of
  $E_K$ or on the $k$-th blow-up $X_k$ of $\mathbb{P}^1$.  In fact we
  will use the projection formula to compute intersections on the
  minimal model $\cE'$ for $E$ after the base change $\Spec
  V[X]/(X^n-\pi)\rightarrow \Spec V$ where $\pi$ is a uniformizing
  parameter of $V$.  A description of the resolution of singularities
  of the base change can be found in \cite[Expos{\'e} 1, Propositio
  2.2]{asterisqueun}.

  On the minimal model $\cE'$, the special fiber has $2kn$ components.
  Let $Z_0$ denote the component of the origin of the elliptic curve,
  and let us denote the other components as $Z_1, \dots, Z_{2kn - 1}$
  in such a way that $Z_i$ meets $Z_{i+1}$ for $0 \leq i \leq (2kn -
  1)$ and $Z_{2k -1}$ meets $Z_0$ (figure 4).

  The divisor of $n$-torsion points meets only the components $Z_i$
  for which $i$ is a multiple of $2k$; the multiplicity of each
  intersection is $n$. The components $Z_j$ for which $j$ is a
  multiple of $n$ are the only ones not contracted by the morphism to
  the plane model $E$.  The contribution at $Q_j$ in the intersection
  number $(C_k.F_j)$ for $j \not= 0, k$ will be
   $$
   n \cdot \big| \{ \text{$m$ such that $(j-1)n \leq 2km \leq jn$}
   \}\big| $$
   Write $n=2kq+r$ with $0 \leq r<2k$. We have
   $$
   \Big| \big| \{ \text{$m$ such that $(j-1)n \leq 2km \leq jn$}
   \}\big| - q \Big| \leq 1.$$
   Thus, we have
   $$(**)\left| (C_k.F_j)- 2n\frac{n-r}{2 k} \right| \leq 2n.$$
   Since
   $(C_k.F_k) = n \frac{n-r}{k}$ we obtain
   $$(H_n.D) \simeq 2 \sum_{j\leq(k-1)}j 2n\frac{n-r}{2k} +
   2k n \frac{n-r}{k}$$
   with an error at most $2\sum_{j\leq(k-1)}j(2n) = 2
   \frac{k(k-1)}{2} 2n$. Hence, we have
$$|(H_n.D) -  (n-r)nk|      \leq k(k-1) 2n,$$
     so
     $$
     \lim_{n\rightarrow\infty} \frac{1}{n^2}(H_n.D)= k.$$
\\[0.4 in]
     \includegraphics[scale=0.6]{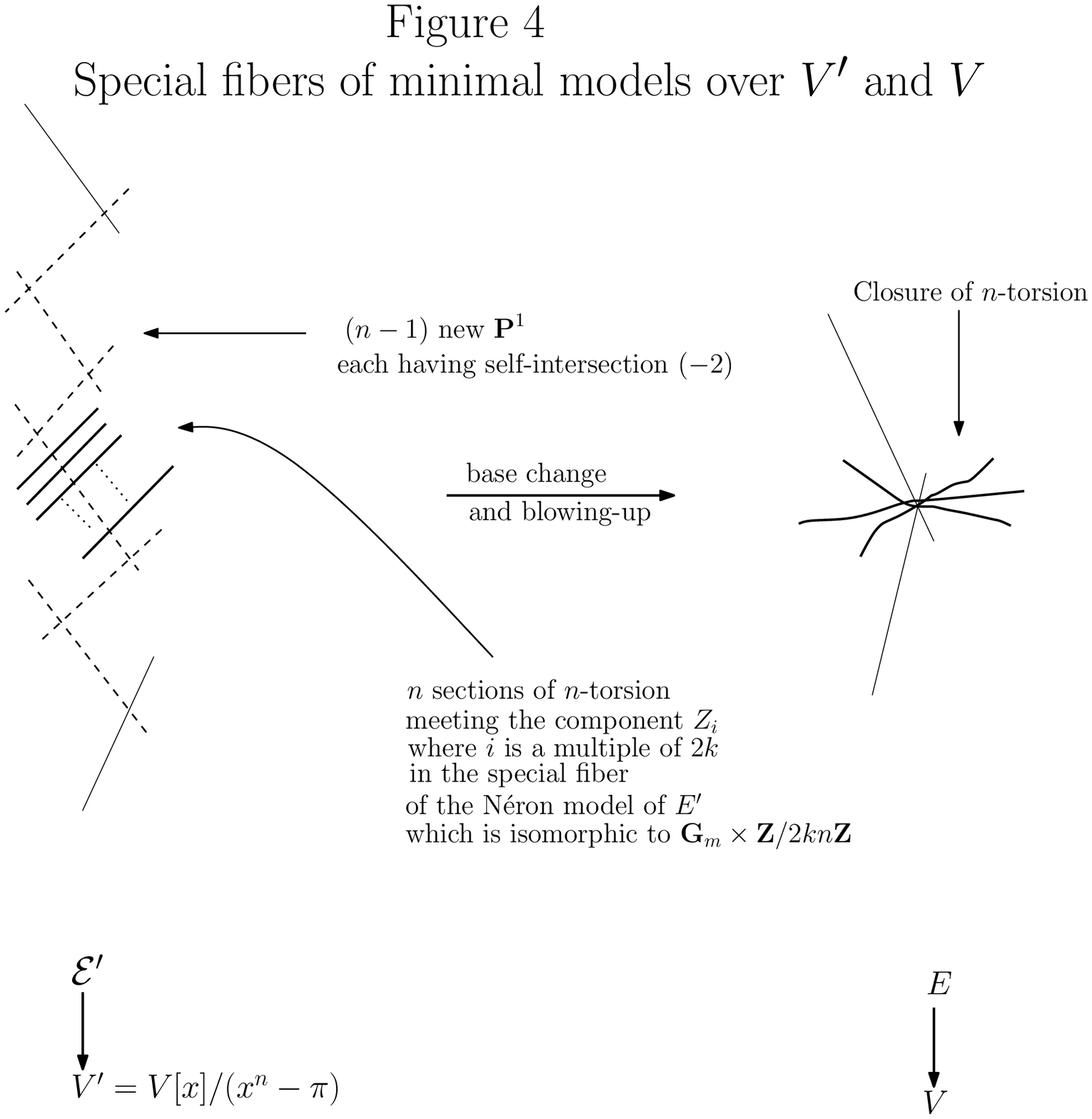}
\\[0.4 in]

     This finishes the proof in the case when $n$ and $2k$ are
     relatively prime.  For the case where $n$ and $2k$ have a gcd $m$
     greater than 1 the formula (*) is still valid. The $n$-torsion
     distribute themselves in packets of $m$ in components of the
     special fiber (see figure 5).  Thus, the estimate (**) for
     $(C_k.F_i)$ has now an error term of at most $m$.
\\[0.4 in]
\includegraphics[scale=0.6]{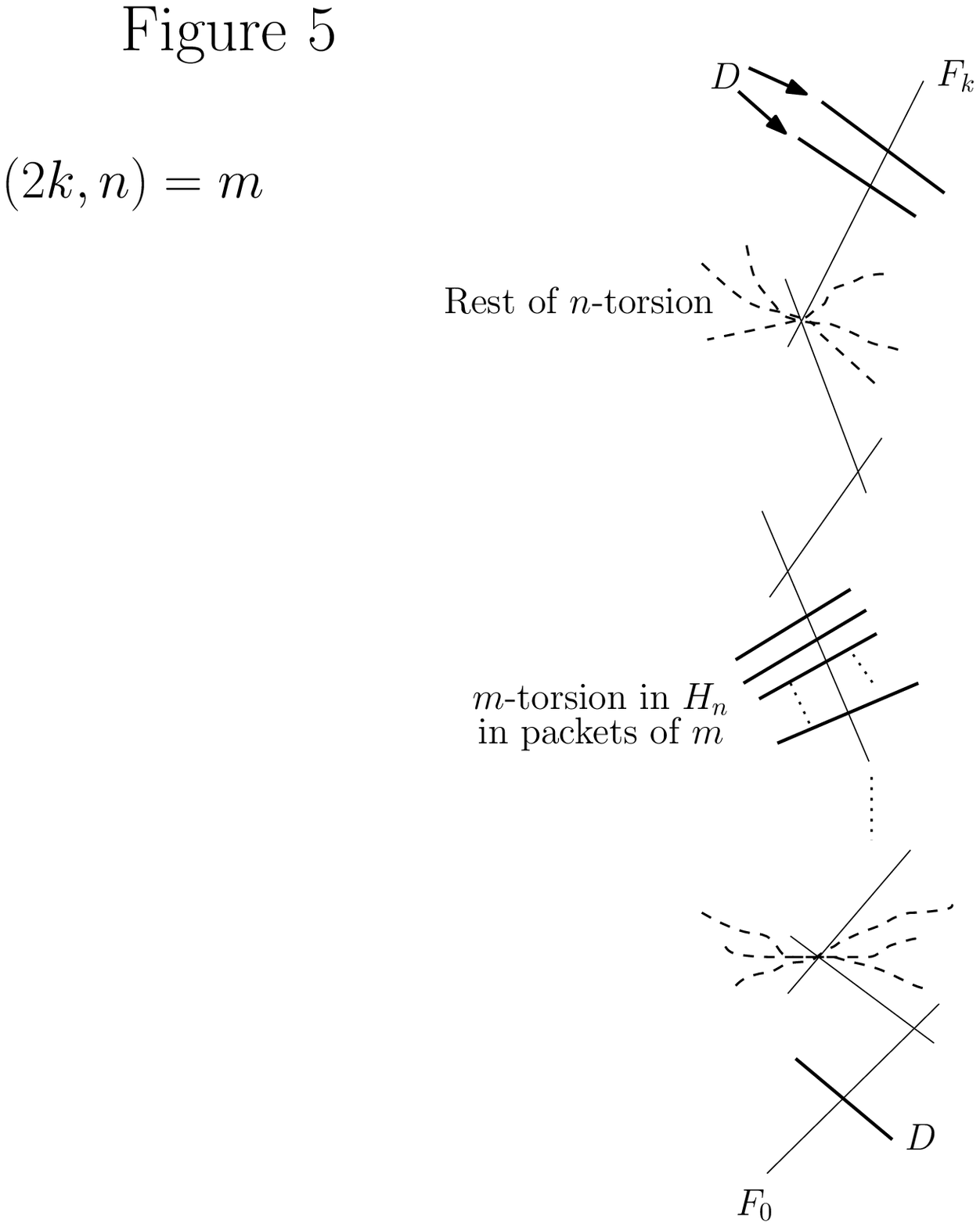}\\
\\[0.4 in]

Adding as before, we now obtain
$$|(H_n.D) - (n-r)nk| \leq \sum_{j\leq(k-1)}2jm = m k(k-1) \leq
k(k-1)(2k).$$
Letting $n$ go to $\infty$ we see again that
$$
\lim _{n\rightarrow \infty}\frac{1}{n^2}(H_n.D)= k.$$
\end{proof}
\vspace{4mm}
\noindent {\it Acknowledgments.}The authors would like to thank
M.~Baker and R.~Rumely for many helpful discussions.

\providecommand{\bysame}{\leavevmode\hbox to3em{\hrulefill}\thinspace}
\providecommand{\MR}{\relax\ifhmode\unskip\space\fi MR }
\providecommand{\MRhref}[2]{%
  \href{http://www.ams.org/mathscinet-getitem?mr=#1}{#2}
}
\providecommand{\href}[2]{#2}

\end{document}